\documentclass[12pt, a4paper]{amsart}
\usepackage{amsmath,amssymb,amscd,amsfonts}

\usepackage{ifthen}
\usepackage[T1]{fontenc}
\usepackage[utf8]{inputenc}
\usepackage[all]{xy}
\usepackage{graphicx}
\usepackage{enumerate}
\usepackage{xspace}
\usepackage{pdfsync}
\usepackage{epic}
\usepackage{dsfont}

\newtheorem{theorem}{Theorem}[section]

\newtheorem{lemma}[theorem]{Lemma}

\def\cO{{\mathcal O}}

\def\cL{{\mathcal L}}

\begin{document}
\title[]{Local projectivity of Lagrangian fibrations on  Hyperk\"ahler manifolds}

\author{Fr\'ed\'eric Campana}
\address{Universit\'e Lorraine \\
 Institut Elie Cartan\\
Nancy \\ Institut Universitaire de France\\ and KIAS scholar,\ KIAS\\
85 Hoegiro, Dongdaemun-gu\\
Seoul 130-722, South Korea,}

\email{frederic.campana@univ-lorraine.fr}

\date{\today}

\maketitle

%\ \ \ \ \ \ \ \ \ \ \ \ \ \ \ \ \ \ \ \ \ \ \ \ \ \ \ \ {\it A Jean-Pierre Demailly}

%\tableofcontents

%%%%%%%%%%%%%%%%%%%%%%%%%%%%%%%%%%%%%%%%%%%%
%%%%%%%%%%%%%%%%%%%%%%%%%%%%%%%%%%%%%%%%%%%%%%

\begin{abstract} We show that if $f:X\to B$ is a Lagrangian fibration from a compact connected K\"ahler hyperk\"ahler manifold $X$ onto a projective normal variety $B$, then $f$ is locally projective. This answers a question raised by L. Kamenova and strengthens a former result (see \cite{Ca05}, Proposition 2.1), according to which the smooth fibres of $f$ are projective.\end{abstract}

%%%%%%%%%%%%%%%%%%%%%%%%%%%%%%%%%%%%%%%%%%%%
%%%%%%%%%%%%%%%%%%%%%%%%%%%%%%%%%%%%%%%%%%%%%%

\section{proof}

The proof consists of simple observations, based on several classical difficult results.

\begin{lemma}\label{locprojcrit} Let $f:X\to B$ be a fibration from a compact connected K\"ahler manifold onto a normal complex projective variety $B$. Assume the existence of a K\"ahler form $w$ on $X$ whose restriction to the generic smooth fibre $X_b$ of $f$ is rational in cohomology (so that $X_b$ is projective, by Kodaira's theorem). Assume additionally that the direct image sheaf $R^2f_*(\cO_X)$ is torsion free. Then $f$ is locally projective.
\end{lemma}

Before giving the proof, let us show how it implies the statement in the abstract. Assume additionally that $X$ is hyperk\"ahler. There thus exists a K\"ahler form $w$ and a holomorphic symplectic form $s$ such that $[w+Re(s)]\in H^2(X,\Bbb Z)$ (the case in which $s$ can be chosen to vanish is trivial, by Kodaira's theorem). This implies the first hypothesis, since the restriction of $w+Re(s)$ to $X_b$ coincides with the restriction of $w$, by the Lagrangian hypothesis. The second hypothesis is fulfilled, by Koll\`ar's theorem (\cite{K}, Theorem 2.1), which asserts that $R^jf_*(K_X)$ is torsionfree, for any $j\geq 0$, and since $K_X=\cO_X$ here, taking $j=2$. The result is stated for $X$ projective in loc.cit., but the proof works in the K\"ahler case as well.

\begin{proof} (of lemma \ref{locprojcrit}): Let $0\in B$ be arbitrary. Choose an open Stein neighborhood of $0\in B$ such that $X_0\subset X_U:=f^{-1}(U)$ is a deformation retract, so that the natural restriction maps $H^{j}(X_U,\Bbb Z)\to H^{j}(X_0,\Bbb Z)$ are isomorphisms (as above, we shall only use the index $j=2$).

We have a natural map derived from the exponential sequence: $$exp:H^2(X_U,\Bbb Z)\to H^2(X_U,\cO_X)=H^0(U,R^2f_*(O_{X_U})),$$ the last equality holds because $U$ is Stein (indeed, by Grauert's coherence theorem, and Theorem B of Cartan-Serre, all arrows in the Leray spectral sequence vanish from $E_2^{p,q}=H^p(U,R^qf_*\cO_X)$ on, since $E_2^{p,q}=0$, for $p>0$, so that $E_r^{p,q}=E_2^{p,q}, \forall r>2$).

Let now $u=[w]_{\vert X_U}$, and $u':=exp(u)\in H^0(U,R^2f_*(\cO_X))$ be the image of $u$. By our first hypothesis, this section vanishes on a Zariski open subset of $U$. By torsionfreeness, it thus vanishes everywhere on $U$, and $u$ is thus the Chern class of some line bundle $\cL_U$ on $U$. We now show that $\cL_U$ is $f$-ample on $X_U$.

Let $Z\subset X_b,b\in U$ be any $d$-dimensional irreducible subvariety, $0\leq d$. We have, for the intersection numbers: $$[w]^d.Z=\cL_{U}^d.Z\geq 1,$$ since $[w]^d.Z>0$, because $w$ is a K\"ahler form, and $[w]$ is integral, which gives a positive lower bound on these intersection numbers. The relative version of Nakai-Moishezon criterion (in the version of  H. Grauert, which does not presuppose algebraicity) thus applies, and concludes the proof.
\end{proof}

{\bf Remark:} After the present text was posted, C. Lehn informed me that he already proved in \cite{L}, Theorem 1.1, the projectivity of each fibre of $f$ in our situation, by a mixed Hodge structure argument.

%\section{Introduction}


\begin{thebibliography}{66}




\bibitem{Ca05} F. Campana. \emph{Isotrivialit\'e de certaines familles K\"ahl\'eriennes de vari\'et\'es non projectives.} Math. Z. 252 (2005), 147-156. \smallskip



\bibitem{K} J. Koll\`ar. \emph{Higher direct images of dualizing sheaves I} Ann. Math. 123 (1986), 11-42 \smallskip

\bibitem{L} C. Lehn. \emph{Deformations of Lagrangian subvarieties of holomorphic symplectic manifolds.} Math. Research Letters 23, 473-497 (2016).


\end{thebibliography}
\end{document}